\documentclass[10pt]{article}
\usepackage{amsfonts}
\usepackage{amsmath}
\usepackage{mathrsfs,multirow}
\usepackage{mathrsfs,amscd,amssymb,amsthm,amsmath,bm,graphicx,psfrag,subfigure}

\makeatletter

\renewcommand{\@seccntformat}[1]{{\csname the#1\endcsname}{\normalsize .}\hspace{.5em}}
\makeatother

\def \[{\begin{equation}}
\def \]{\end{equation}}

\newtheorem{thm}{Theorem}[section]

\newtheorem{lem}[thm]{Lemma}
\newtheorem{cor}[thm]{Corollary}

\newenvironment{wst}
{\setlength{\leftmargini}{1.5\parindent}
 \begin{itemize}
 \setlength{\itemsep}{-1.1mm}}
{\end{itemize}}

\begin{document}
\setlength{\baselineskip}{15pt}
\begin{center}{\Large \bf On the positive and negative inertia of weighted graphs\footnote{Financially supported by the National Natural
Science Foundation of China (Grant Nos. 11071096, 11271149).}}

\vspace{4mm}

{\large Shuchao Li\footnote{E-mail: lscmath@mail.ccnu.edu.cn (S.C.
Li),\, 928046810@qq.com (F.F. Song)},\ Feifei Song}\vspace{2mm}

Faculty of Mathematics and Statistics,  Central China Normal
University, Wuhan 430079, P.R. China
\end{center}
\noindent {\bf Abstract}: The number of the positive, negative and zero eigenvalues in the spectrum of the (edge)-weighted graph $G$ are called positive inertia index, negative inertia index and nullity of the weighted graph $G$, and denoted by $i_+(G)$, $i_-(G)$, $i_0(G)$, respectively. In this paper, the positive and negative inertia index of weighted trees, weighted unicyclic graphs and weighted bicyclic graphs are discussed, the methods of calculating them are obtained.

\vspace{2mm} \noindent{\it Keywords}: Positive inertia index; Negative inertia index; Weighted graph; Tree; Unicyclic graph; Bicyclic graph

\vspace{2mm}

\noindent{AMS subject classification:} 05C50, \ 15A18

 {\setcounter{section}{0}
\section{\normalsize Introduction}\setcounter{equation}{0}
In this paper, we only consider simple weighted graphs on positive weight sets. Let $G$ be a weighted graph with vertex set $\{v_1, v_2, \ldots , v_n\}$, edge set $E_G\neq \emptyset$ and $W(G)=\{w_{j}> 0, j=1, 2, \ldots , |E_G|\}$. The function $w: E_G \rightarrow W(G)$ is called a \textit{weight function} of $G$. It is obvious that each weighted graph corresponds to a weight function. Such weighted graph is usually denoted by $G=(V_G,E_G,w)$. If there is no fear of confusion, we often simply call $G$ a weighted graph.

The weighted adjacency matrix $A_w(G)$ of a weighted graph $G$ of order $n$ is defined as the matrix $A(G) = A_w(G) = (a_{ij})$ of order $n$ with $a_{ij}=w(v_iv_j)$ if  $v_iv_j\in E_G$ and $0$ otherwise.
In particular, if each nonzero $a_{ij}=1$, $A(G)$ is the adjacency matrix of its underlying graph of $G$. The \textit{inertia} of $G$ is defined to be the triplet In$(G)=(i_+(G),\,i_-(G),\, i_0(G))$, where $i_+(G)$,\, $i_-(G)$, \,$i_0(G)$ are the number of the positive, negative and zero eigenvalues of $A(G)$ including multiplicities, respectively. $i_+(G)$,\, and $i_-(G)$ are called the \textit{positive, negative index} of inertia (for short, positive, negative index) of $G$, respectively. Traditionally, $i_0(G)$ is called the \textit{nullity} of $G$. Obviously $i_+(G)+i_-(G)+ i_0(G)=n$. If $H$ is a subgraph of $G$ with $w_{H}(e)=w_{G}(e)$ for each $e\in E_H$, then $H$ is called a \textit{weighted subgraph} of $G$. For an induced weighted subgraph $H$ of a weighted graph $G$, let $G-H$ be the subgraph obtained from $G$ by deleting all vertices of $H$ and all incident edges. We define the union of two disjoint weighted graphs $G_1=(V_1,E_1,w_1)$ and $G_2=(V_2,E_2,w_2)$, denoted by $G_1\cup G_2$, as the graph with vertex set $V_1\cup V_2$, edge set $E_1\cup E_2$ and weight function $w:E_1\cup E_2\rightarrow W(G_1)\cup W(G_2)$ satisfying
$$
w(e)  =\left\{
         \begin{array}{ll}
          w_1(e), & \hbox{if $e\in E_1;$} \\
           w_2(e), & \hbox{if $e\in E_2;$} \\
           0, & \hbox{otherwise.}
         \end{array}
       \right.
$$

For a subgraph $H$ of $G$ and a vertex $x \in V_G$, we denote $N_{H}(x)= N_{G}(x)\bigcap V_H$. A vertex of a graph $G$ is called \textit{pendant} if it has degree one. A tree (resp. unicyclic graph, bicyclic graph) is a simple connected graph in which the number of edges equals the number of vertices plus $-1\,$ (resp. 0, 1). A complete graph, a path, a cycle and a star of order $n$ are denoted by $K_n, P_n, C_n$, and $K_{1,n-1}$, respectively. An isolated vertex is sometime denoted by $K_1$.

The characteristic polynomial of $A(G)$ is said to be characteristic polynomial of $G$, denoted by $\phi(G,\lambda)$ (or, $\phi(G)$ for short). The coefficient $a_0(G)$, $a_1(G)$ are the constant term, one degree term of $\phi (G,\lambda)$, respectively.

The inertia of unweighted graphs has attracted some attention. Gregory et al. \cite{18} studied
the subadditivity of the positive, negative indices of inertia and developed certain properties of
Hermitian rank which were used to characterize the biclique decomposition number. Gregory
et al. \cite{19} investigated the inertia of a partial join of two graphs and established a few relations
between inertia and biclique decompositions of partial joins of graphs. Daugherty \cite{28} characterized
the inertia of unicyclic graphs in terms of matching number and obtained a linear-time
algorithm for computing it. Ma et al. \cite{80} studied the methods of calculating the positive and negative inertia index of unweighted graphs. Yu at al. \cite{17} investigated the minimal positive index of inertia
among all unweighted bicyclic graphs of order $n$ with pendant vertices, and characterized the
bicyclic graphs with positive index 1 or 2. The nullity of
unweighted graphs has been studied well in the literature, one may be referred to \cite{11,9,8} and the survey \cite{20}. There is also a
large body of knowledge regarding to the inertia of unweighted graphs
due to its many applications in chemistry (see \cite{3,21,22,16} for details). The study of eigenvalues of weighted graph
attracts much attention. Several results about the (Laplacian) spectral radius of weighted graphs were derived, one may be referred to \cite{23,24,27,26,25}.
Only a few papers considered the inertia of weighted graphs. Yu et al. determined a lower bound on the positive, negative index of weighted $n$-vertex unicyclic graphs with fixed girth and characterize all weighted unicyclic graphs attaining this lower bound. Motivated by the above description, we shall focus on the method on computing the inertia number of weighted trees, weighted unicyclic graphs and weighted bicyclic graphs of order $n$, respectively.

This paper is organized as follows: in Section 2, some necessary lemmas are given. In Section
3, the methods of calculating the positive and negative inertia index of weighted trees and weighted unicyclic graphs are obtained. In the last section, the methods of calculating the positive and negative inertia index of weighted bicyclic graphs are obtained.

\section{\normalsize  Preliminary results}\setcounter{equation}{0}
In this section, we cite some previous results. Suppose $A$, $B$ are two Hermitian matrices of order $n$, if there exists an invertible matrix $P$ of order $n$ such that $P^{\ast}AP=B$, $P^{\ast}$ denotes the conjugate transpose of $P$, then we say that $A$ is \textit{congruent} to $B$, denoted by $A\cong B$.

\begin{lem}[\cite{9}]\label{lem2.1}
Let $A$, $B$ be two Hermitian matrices of order n, such that $A\cong B$. Then $i_+(A)=i_+(B)$, $i_-(A)=i_-(B)$, $i_0(A)=i_0(B)$.
\end{lem}

It is easy to obtain the following result.
\begin{lem}\label{lem2.2}
Let $G=G_1\cup G_2\cup \ldots \cup G_t$ be a weighted graph, where $G_i$\, $(i=1,2,\ldots, t)$ are connected components of $G$. Then $i_+(G)=\bigcup\limits_{i=1}^{t}i_+(G_i)$, $i_-(G)=\bigcup\limits_{i=1}^{t}i_-(G_i)$, $i_0(G)=\bigcup\limits_{i=1}^{t}i_0(G_i)$.
\end{lem}

Let $M$ be a Hermitian matrix. We denoted three types of elementary congruence matrix operations (ECMOs) on $M$ as follows:
\begin{enumerate}
                                                           \item  interchanging $i$th and $j$th rows of $M$, while interchanging $i$th and $j$th columns of $M$;                 \item multiplying $i$th row of $M$ by a non-zero number $k$, while multiplying $i$th column of $M$ by $k$;
                                                           \item adding $i$th row of $M$ multiplied by a non-zero number $k$ to $j$th row, while adding $i$th column of $M$ multiplied by $k$ to $j$th column.
\end{enumerate}
By Lemma $2.1$, the ECMOs do not change the inertia of a Hermitian matrix.

By interlacing inequalities for eigenvalues of Hermitian matrices, we can deduce the following result.
\begin{lem}[\cite{17}]\label{lem2.3}
Let $A$ be an $n\times n$ Hermitian matrix and $B$ be the Hermitian matrix obtained by bordering $A$ as followings:
 \[B=
\left(\begin{array}{cccc}
A & y\\
y^{\ast} & a \notag
\end{array}\right),
\]
where $y$ is a column vector, $y^{\ast}$ denotes the conjugate transpose of $y$ and $a$ is a real number. Then $i_+(B)-1 \leq i_+(A)\leq i_+(B)$, $i_-(B)-1 \leq i_-(A)\leq i_-(B)$.
\end{lem}

 The following result is a direct consequence of Lemma 2.3.
\begin{lem}[\cite{09}]\label{lem2.4}
Let $H$ be an induced weighted subgraph of $G$. Then $i_+(H) \leq i_+(G)$, $i_-(H) \leq i_-(G)$.
\end{lem}

\begin{lem}[\cite{09}]\label{lem2.5}
Let $G=(V_G,E_G,w)$ be a weighted graph containing a pendant vertex $v$ with its unique neighbor $u$. Then $i_+(G)=i_+(G-u-v)+1$, $i_-(G)=i_-(G-u-v)+1$.
\end{lem}
\begin{proof}
We present the proof here for completeness. Assume that all vertices in $V_G$ are indexed by $\{v_1,v_2,\ldots ,v_n\}$ with $v_1=v$, $v_2=u$, $w(uv)=a> 0$. Then it follows that
 \[A(G)=
\left(\begin{array}{cccccc}
0 & a& 0 & \ldots & 0\\
a & 0& a_{23} & \ldots & a_{2n} \\
0 & a_{32}& 0 & \ldots & a_{3n} \\
\vdots & \vdots & \vdots & \ddots & \vdots\\
0 & a_{n2}& a_{n3} & \ldots & 0 \notag
\end{array}\right),\]
where the first two rows and columns are labeled by $v_1$, $v_2,$ respectively.

By applying the ECMOs on $A(G)$, it is easy to show that $A(G)$ is congruent to
 \[A(P_2\cup(G-u-v))=
\left(\begin{array}{cccccc}
0 & a& 0 & \ldots & 0\\
a & 0& 0 & \ldots & 0 \\
0 & 0& 0 & \ldots & a_{3n} \\
\vdots & \vdots & \vdots & \ddots & \vdots\\
0 & 0& a_{n3} & \ldots & 0 \notag
\end{array}\right).\]

Therefore $i_+(G)=i_+(G-u-v)+i_+(P_2)=i_+(G-u-v)+1$, $i_-(G)=i_-(G-u-v)+i_-(P_2)=i_-(G-u-v)+1$.
\end{proof}
\begin{lem}\label{lem2.6}
Let $G=(V_G,E_G,w)$ be a weighted graph containing path with four vertices of degree 2 as depicted in Fig. 1, where the numbers on the edges denote the weights of the corresponding edges. Let $H$ be obtained from $G$ by replacing this path with an edge, whose weight is $\frac{a_1a_3a_5}{a_2a_4}$. Then $i_+(G)=i_+(H)+2$, $i_-(G)=i_-(H)+2$.
\end{lem}
\begin{figure}[h!]
\begin{center}
  \psfrag{1}{$v_6$}\psfrag{2}{$v_3$}\psfrag{3}{$v_2$}\psfrag{4}{$v_1$}
  \psfrag{5}{$v_4$}\psfrag{6}{$v_5$}\psfrag{a}{$a_1$}\psfrag{b}{$a_2$}
  \psfrag{c}{$a_3$}\psfrag{d}{$a_4$}\psfrag{e}{$a_5$}\psfrag{f}{$\frac{a_1a_3a_5}{a_2a_4}$}
  \psfrag{7}{$v_6$}\psfrag{8}{$v_5$}\psfrag{g}{$G$}\psfrag{h}{$H$}
  \includegraphics[width=100mm]{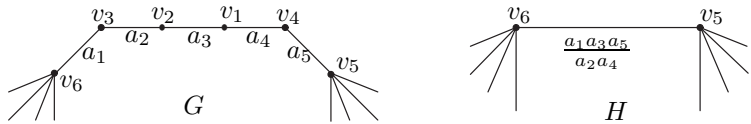}\\
  \caption{The weighted graph operation I of calculating positive and negative inertia index.}
\end{center}
\end{figure}
\begin{proof}
The adjacency matrix $A(G)$ has the following form
 \[A(G)=
\left(\begin{array}{ccccccccc}
0 & a_3& 0 & a_4& 0&0&0&\ldots & 0\\
a_3 & 0& a_2 &0&0&0& 0&\ldots & 0 \\
0 & a_2& 0 & 0& 0&a_1&0&\ldots & 0\\
a_4 & 0& 0 &0&a_5&0& 0&\ldots & 0 \\
0 & 0& 0 &a_5&0&0& a_{57}&\ldots & a_{5n} \\
0 & 0& a_1 & 0& 0&0&a_{67}&\ldots & a_{6n}\\
0 & 0& 0 &0&a_{75}&a_{76}& 0&\ldots & a_{7n} \\
\vdots & \vdots& \vdots & \vdots& \vdots&\vdots&\vdots&\ddots& \vdots\\
0& 0 &0&0&a_{n5}& a_{n6}&a_{n7}&\ldots & a_{nn} \notag
\end{array}\right).\]

It is routine to check that $A(G)$ is congruent to
 \[A(P_2\cup P_2\cup H)=
\left(\begin{array}{ccccccccc}
0 & a_3& 0 & 0& 0&0&0&\ldots & 0\\
a_3 & 0& 0 &0&0&0& 0&\ldots & 0 \\
0 & 0& 0 &\frac{-a_2a_4}{a_3} & 0&0&0&\ldots & 0\\
0 & 0& \frac{-a_2a_4}{a_3}&0&0&0& 0&\ldots & 0 \\
0 & 0& 0 &0&0&\frac{a_1a_3a_5}{a_2a_4}& a_{57}&\ldots & a_{5n} \\
0 & 0& 0 & 0& \frac{a_1a_3a_5}{a_2a_4}&0&a_{67}&\ldots & a_{6n}\\
0 & 0& 0 &0&a_{75}&a_{76}& 0&\ldots & a_{7n} \\
\vdots & \vdots& \vdots & \vdots& \vdots&\vdots&\vdots&\ddots& \vdots\\
0& 0 &0&0&a_{n5}& a_{n6}&a_{n7}&\ldots & a_{nn} \notag
\end{array}\right).\]

Therefore $i_+(G)=i_+(H)+2$, $i_-(G)=i_-(H)+2$.
\end{proof}

\section{\normalsize Positive and negative inertia index of weighted trees and weighted unicyclic graphs}
An edge subset $M\subseteq E_G$ is called a \textit{matching} of $G$ if no two edges of $M$ share a common vertex. A matching $M$ is called maximum in $G$ if it has maximum cardinality among all matchings of $G$, and is called perfect if every vertex of $G$ is incident with (exactly) one edge in $M$. Obviously, a perfect matching is a maximum matching. The cardinality of a maximum matching is called the \textit{matching number} of $G$, denoted by $q(G)$.
\begin{lem}\label{lem3.1}
Let $T$ be a weighted tree of order $n$. Then $i_+(T)=i_-(T)=q(T)$.
\end{lem}
\begin{proof}
We apply induction on $n$. If $n=1$, $i_+(T)=i_-(T)=q(T)=0$. So suppose the assertion holds for smaller values of $n$. Assume $T$ be a weighted tree containing a pendant vertex $v$ with its unique neighbor $u$. Then by Lemma 2.5, we have
$$
 i_+(T)=i_+(T-u-v)+1,\ \  i_-(T)=i_-(T-u-v)+1.
$$
Assume $N_{T-v}(u)=\{u_1,u_2,\ldots,u_t\}$ $(t\geq 1)$, and $T_1,T,_2,\ldots ,T_t$ are the components of $T-u-v$ that contain the vertices $u_1,u_2,\ldots,u_t$, respectively. So by Lemma 2.2 and induction,
\begin{eqnarray*}
   && i_+(T)=i_+(T-u-v)+1=\sum_{i=1}^{t}i_+(T_i)+1=\sum_{i=1}^{t}q(T_i)+1=q(T-u-v)+1=q(T),\\
  && i_-(T)=i_-(T-u-v)+1=\sum_{i=1}^{t}i_-(T_i)+1=\sum_{i=1}^{t}q(T_i)+1=q(T-u-v)+1=q(T).
\end{eqnarray*}

This completes the proof.
\end{proof}

By Lemmas \ref{lem2.2} and \ref{lem3.1} we have the following results immediately.
\begin{cor}
Let $G$ be a weighted acyclic graph (or a forest) of order $n$. Then $i_+(G)=i_-(G)=q(G)$.
\end{cor}

A tree is called a PM-tree if it has a perfect matching.
\begin{cor}
Let $T$ be a weighted PM-tree of order $n$. Then $i_+(T)=i_-(T)=\frac{n}{2}$.
\end{cor}
\begin{lem}\label{lem3.2}
Let $C_n=v_1v_2\ldots v_nv_1$ be a weighted cycle of order $n$, $w(v_iv_{i+1})=a_i$ $(1\leq i\leq n)$ and let $v_{n+1}=v_1$. If $n\equiv 0\pmod{4}$, then
$$
i_+ (C_n)= i_- (C_n)= \left\{
  \begin{array}{ll}
    \frac{n}{2}-1, & \hbox{if $\prod_{i=1}^{\frac{n}{2}}a_{2i-1}= \prod_{i=1}^{\frac{n}{2}}a_{2i};$} \\[4pt]
    \frac{n}{2}, & \hbox{otherwise.}
  \end{array}
\right.
$$
\end{lem}
\begin{proof}
Suppose $C_4^*$ is a weighted cycle of order 4 and all vertices in $C_4^*$ are indexed by $\{v_1,v_{n-2},v_{n-1},v_n\}$ with
$$
w(v_{n-2}v_{n-1})=a_{n-2}, \, w(v_{n-1}v_n)=a_{n-1},\, w(v_{n-1}v_n)=a_n,\, w(v_1v_{n-2})= \frac{ \prod_{i=1}^{ \frac{n}{2}-1}a_{2i-1}}{ \prod_{i=1}^{\frac{n}{2}-2}a_{2i}}.
$$
Then applying Lemma \ref{lem2.6} repeatedly, we have $i_+(C_n)= 2(\frac{n}{4}-1)+i_+(C_4^{\ast})$ and $i_-(C_n)= 2(\frac{n}{4}-1)+i_-(C_4^{\ast})$.

The adjacency matrix $A(C_4^*)$ has the following from
\[
A(C_4^{\ast})=
\left(\begin{array}{cccc}
0 & \frac{\prod_{i=1}^{\frac{n}{2}-1}a_{2i-1}}{\prod_{i=1}^{\frac{n}{2}-2}a_{2i}}& 0  & a_n\\[3pt]
\frac{\prod_{i=1}^{\frac{n}{2}-1}a_{2i-1}}{\prod_{i=1}^{\frac{n}{2}-2}a_{2i}} & 0& a_{n-2}  & 0 \\[3pt]
0 & a_{n-2}& 0 &  a_{n-1} \\[3pt]
a_n & 0& a_{n-1} & 0 \notag
\end{array}\right).
\]
Then we have
\begin{eqnarray*}
  a_0(C_4^*) &=& \left(a_na_{n-2}- a_{n-1}\frac{\prod_{i=1}^{\frac{n}{2}-1}a_{2i-1}}{\prod_{i=1}^{\frac{n}{2}-2}a_{2i}}\right)^{2}\geq 0, \\
  a_1(C_4^{\ast}) &=& 0.
\end{eqnarray*}

Note that $C_4^{\ast}$ contains $P_2$ as an induced subgraph, hence $i_+(C_4^{\ast})\geq 1$ and $i_-(C_4^{\ast})\geq 1.$ Therefore
$$
i_+ (C_4^{\ast})= i_- (C_4^{\ast})=\left\{
  \begin{array}{ll}
  1  , & \hbox{if $\prod_{i=1}^{\frac{n}{2}}a_{2i-1}= \prod_{i=1}^{\frac{n}{2}}a_{2i}$;} \\[5pt]
  2  , & \hbox{otherwise.}
  \end{array}
\right.
$$

This completes the proof.
\end{proof}
\begin{lem}\label{lem3.3}
Let $C_n$ be a weighted cycle of order $n$. Then
\[
i_+ (C_n)= \left\{ {\begin{array}{*{20}c}
   \frac{n+1}{2},  & \hbox{if $n\equiv 1\pmod{4}$};  \\[5pt]
   \frac{n}{2}, & \hbox{if $n\equiv 2\pmod{4}$}; \notag \\[5pt]
   \frac{n-1}{2}, & \hbox{if $n\equiv 3\pmod{4}$}.  \notag
\end{array}} \right.
i_- (C_n)= \left\{ {\begin{array}{*{20}c}
   \frac{n-1}{2},  & \hbox{if $n\equiv 1\pmod{4}$};  \\[5pt]
   \frac{n}{2}, & \hbox{if $n\equiv 2\pmod{4}$};  \notag \\[5pt]
   \frac{n+1}{2}, & \hbox{if $n\equiv 3\pmod{4}$}.  \notag
\end{array}} \right.
\]
\end{lem}
\begin{proof}
 First we shall consider a weighted cycle $C_3$ with $W(C_3)= \{a_1,a_2,a_3\}$. Then $a_0(C_3)=-2a_1a_2a_3<0$. And note that $i_+(C_3)\geq 1$ and $i_-(C_3)\geq 1$ since $C_3$ contains $P_2$ as an induced subgraph. Therefore $i_+(C_3)=1$, $i_-(C_3)=2$. By the similar argument, we can obtain that $i_+(C_5)= 3$, $i_-(C_5)= 2$, $i_+(C_6)=i_-(C_6)=3$.

 When $n\geq 7$ and $ n\equiv 1,2,3 \pmod{4}$, we can easily obtain the conclusions by Lemma \ref{lem2.6}.
\end{proof}

For a tree $T$ on at least two vertices, a vertex $v\in T$ is called \textit{mismatched} in $T$ if there exists a maximum matching $M$ of $T$ that does not cover $v$; otherwise, v is called \textit{matched} in $T$. If a tree consists
of only one vertex, then this vertex is considered mismatched.

A vertex of a graph $G$ is called \textit{quasi-pendant} if it is adjacent to a pendant vertex. We have the following properties of matched or mismatched vertices of a weighted tree.

\begin{lem}[\cite {11}]\label{lem3.4} Given a tree $T$,
\begin{wst}
\item[{\rm (i)}]if $v$ is a quasi-pendant vertex of $T$, then $v$ is matched in $T;$
\item[{\rm (ii)}]if $v$ is a mismatched vertex of $T$, then for any neighbor $u$ of $v$, $u$ is matched in $T$, and it is also
matched in the component of $T-v$ that contains $u$.
\end{wst}
\end{lem}
\begin{lem}\label{lem3.5}
Let $T$ be a weighted tree with $v\in V_T$. The following statements are equivalent:
\begin{wst}
\item[{\rm (i)}]$v$ is mismatched in $T$.
\item[{\rm (ii)}]$q(T - v) = q(T)$.
\item[{\rm (iii)}]$i_+(T - v) = i_+(T)$
\item[{\rm (iv)}]$i_-(T - v) = i_-(T)$.
\end{wst}
\end{lem}
\begin{proof}
 (i), (ii) have been proved in [Lemma 2.6, 12]. Note $T$ and $T - v$ are
both weighted acyclic graphs, then we know that (iii) and (iv) are obtained by Corollary 3.2.
\end{proof}
Let $G_1$ be a graph containing a vertex $u$, and let $G_2$ be a graph of order $n$ that is disjoint from $G_1$. For
$1 \leq k \leq n$, the $k$-joining graph of $G_1$ and $G_2$ with respect to $u$ is obtained from $G_1 \bigcup G_2 $ by joining $u$
and any $k$ vertices of $G_2$; we denote it by $G_1(u)\odot^kG_2$. Note that the graph $G_1(u)\odot ^kG_2$ is not uniquely
determined when $n > k$. If $G_1= (V_1,E_1,w_1)$ and $G_2= (V_2,E_2,w_2)$, that is, $G_1$ and $G_2$ are weighted graphs, assume $N_{G_2}(u)=\{v_1,v_2,\ldots ,v_k\}$. Then the weight function of $G_1(u)\odot ^kG_2$ is $w: E_{G_1(u)\odot ^kG_2}\rightarrow W(G_1(u)\odot ^kG_2)$ with
$$
w(e)  =\left\{
                                        \begin{array}{ll}
                                          w_1(e), & \hbox{if $e\in E_1$;} \\
                                          w_2(e), & \hbox{if $e\in E_2$;} \\
                                          c_i, & \hbox{if $e=uv_i\,(1\leq i\leq k),$ \,where $c_i>0;$} \\
                                          0, & \hbox{otherwise.}
                                        \end{array}
                                      \right.
$$

\begin{thm}\label{thm3.6}
Let $T$ be a weighted tree with a matched vertex $u$ and let $G$ be a weighted graph of order $n$. Then for each positive
integer $k$ $(1 \leq k \leq n)$,
$i_+(T(u)\odot ^kG) = i_+(T) + i_+(G),\ \ i_-(T(u) \odot ^kG) = i_-(T) + i_-(G)$.
\end{thm}
\begin{proof}
Now we prove the $i_+(T(u)\odot ^kG) = i_+(T) + i_+(G)$ by induction on the matching number $q(T)$.

If $q(T) = 1$, then $T$ is a star graph $K_{1,s+1}$ $(s \geq 0)$, and $u$ is the unique quasi-pendant vertex of $T$.
Suppose $v$ is a pendant of $T$ that is adjacent to $u$. By Lemmas 2.2 and 2.5, we have
$$
i_+(T(u)\odot ^kG) = i_+\left((T(u)\odot ^kG)-u-v\right)+1=i_+(sK_1\cup G) +1= i_+(T) + i_+(G).
$$

Suppose $i_+(T(u)\odot ^kG) = i_+(T) + i_+(G)$ holds for any weighted tree $T$ with $q(T) \leq t$\, $(t \geq 1)$. Now we consider a tree $T$ with
$q(T) = t + 1 \geq 2$. As $q(T) \geq 2$, we know that $T$ contains a pendant vertex $v$ and its unique neighbor
$w$, where $v, w$ are both different from $u$. Let $T_1 = T-v-w$. Then $q(T_1) = t$ and
$u$ is matched in $T_1$. By Lemma 2.5 and induction, we have
\begin{center}
$i_+(T(u)\odot ^kG) = i_+((T(u)\odot ^kG)-w-v)+1=i_+(T_1(u)\odot ^kG) +1= i_+(T_1) + i_+(G)+1=i_+(T) + i_+(G)$.
\end{center}

By the similar argument, we can obtain $i_-(T(u) \odot ^kG) = i_-(T) + i_-(G)$, which is omitted here.
\end{proof}
\begin{thm}\label{thm3.7}
Let $T$ be a weighted tree with a mismatched vertex $u$ and let $G$ be a weighted graph of order $n$. Then for each positive
integer $k$ $(1 \leq k \leq n)$, one has
\begin{eqnarray*}
  i_+(T(u)\odot ^kG) &=& i_+(T-u) + i_+(G+u)= i_+(T) + i_+(G+u), \\
  i_-(T(u)\odot ^kG) &=& i_-(T-u) + i_-(G+u)= i_-(T) + i_-(G+u),
\end{eqnarray*}
where $G+u$ is the induced weighted subgraph of $T(u)\odot ^kG.$
\end{thm}
\begin{proof}
In the tree
$T$, suppose that $u_1, u_2, \ldots , u_{m}$ $(m \geq 1)$ are all neighbors of $u$, and $T_1, T_2, \ldots , T_{m}$ are the components
of $T - u$ that contain the vertices $u_1, u_2, \ldots , u_{m}$, respectively. By Lemma 3.6(ii), we know that $u_i$ $(i =
1, 2, \ldots ,m)$ are both matched in $T$ and $T_i$ . Applying Theorem 3.8 repeatedly, we have
\begin{eqnarray*}
  i_+(T(u)\odot ^kG) &=& \sum_{i=1}^{m}i_+(T_i) + i_+(G+u)=i_+(T-u) + i_+(G+u)= i_+(T) + i_+(G+u), \\
  i_-(T(u)\odot ^kG) &=& \sum_{i=1}^{m}i_-(T_i) + i_-(G+u)=i_-(T-u) + i_-(G+u)= i_-(T) + i_-(G+u),
\end{eqnarray*}
as desired.
\end{proof}

In the end of this section, we deal with weighted unicyclic graphs. Let $G$ be a weighted unicyclic graph and let $C$ be the
unique cycle of $G$. For each vertex $v \in C$, let $G\{v\}$ be the induced weighted connected subgraph of $G$ with the
maximum possible number of vertices, which contains the vertex $v$ and contains no other vertices of
$C$. Then for all vertices $v \in C$, $G\{v\}$ is a weighted tree and $G$ is obtained by identifying the vertex $v$ of $G\{v\}$ with
the vertex $v$ on $C$. The unicyclic graph $G$ is called of Type I, if there exists a vertex $v$ on the cycle such
that $v$ is matched in $G\{v\}$; otherwise, $G$ is called of Type I\!I.

If $G$ is of Type I, then $G = G\{v\}(v) \odot^{2} (G- G\{v\})$ for some matched vertex $v$ of $G\{v\}$, where
$G\{v\}$ and $G - G\{v\}$ are both nontrivial weighted trees. Thus we know $i_+(G) = i_+(G\{v\}) + i_+(G - G\{v\})$, $i_-(G) =
i_-(G\{v\}) + i_-(G - G\{v\})$ by Theorem 3.8.
If $G$ is of Type I\!I and $G $ is not a cycle, suppose $G\{v\}$ contains vertices other than $v$, every neighbor of
$v$ in $G\{v\}$ is matched in the component of $ G\{v\}- v$ that contains the neighbor. By Theorem 3.9, we have
\begin{eqnarray*}
  &&i_+(G) = i_+(G\{v\} - v) + i_+(G - G\{v\} + v), \ \  i_-(G) = i_-(G\{v\}- v) + i_-(G - G\{v\} + v).
\end{eqnarray*}
Applying Theorem 3.9 repeatedly, we have
\begin{eqnarray*}
  i_+(G) &=& \sum \limits_{v\in C}i_+(G\{v\} - v) + i_+(C)= i_+(G - C) + i_+(C), \\
  i_-(G) &=& \sum \limits_{v\in C}i_-(G\{v\} - v) + i_-(C)= i_-(G- C) + i_-(C).
\end{eqnarray*}

Hence, we have
\begin{thm}\label{thm3.15}
Let $G$ be a weighted unicyclic graph and let $C$ be the cycle of $G$.
\begin{wst}
\item[{\rm (i)}]
If $G$ is of Type I and let $v \in C$ is matched
in $G\{v\}$, then
$$
i_+(G) =i_+(G\{v\}) + i_+(G - G\{v\}), \ i_-(G) =i_-(G\{v\}) + i_-(G - G\{v\}).
$$
\item[{\rm (ii)}]
If $G$ is of Type II, then
$$
i_+(G) = i_+(G - C) + i_+(C),\ i_-(G) = i_-(G- C) + i_-(C).
$$
\end{wst}
\end{thm}
\section{\normalsize  Positive and negative inertia index of weighted bicyclic graphs}\setcounter{equation}{0}
Let $C_{p}=u_1u_2\ldots u_{p}u_1$ and $C_{q}=v_1v_2\ldots v_{q}v_1$ be two vertex-disjoint cycles. A $\infty $-graph $\infty (p, l, q)$ is
obtained from $C_{p}$ and $C_{q}$ by joining $u_1$ and $v_1$ by a path $v_1w_2 \ldots w_{l-1} v_1$ of length $l-1$, where $l\geq 1$; $l = 1$
means identifying $u_1$ with $v_1$. Let $P_{p}, P_{l}$, and $P_{q}$ be three vertex-disjoint paths, where $\min\{p, l, q\} \geq 2$
and at most one of them is 2. Identifying the three initial vertices and terminal vertices of them,
respectively, the resultant graph is called a $\theta$-graph, denoted by $\theta (p, l, q)$. The weighted  $\infty(p, l, q)$ and $\theta(p, l, q)$ are depicted in Fig. 2, where the numbers on the edges denote the weights of the corresponding edges.
\begin{figure}[h!]
\begin{center}
  \psfrag{1}{$u_1$}\psfrag{3}{$u_3$}\psfrag{d}{$u_{i+1}$}\psfrag{c}{$u_i$}
  \psfrag{e}{$u_{p}$}\psfrag{a}{$v_{q}$}\psfrag{2}{$w_2$}\psfrag{3}{$w_3$}\psfrag{4}{$w_{l-2}$}\psfrag{5}{$w_{l-1}$}
  \psfrag{T}{$a_1$}\psfrag{V}{$a_i$}\psfrag{U}{$a_{p}$}\psfrag{Q}{$b_1$}
  \psfrag{S}{$b_i$}\psfrag{R}{$b_{q}$}\psfrag{8}{$v_i$}\psfrag{6}{$v_1$}
  \psfrag{7}{$v_2$}\psfrag{9}{$v_{i+1}$}\psfrag{M}{$c_1$}\psfrag{N}{$c_2$}
  \psfrag{O}{$c_{l-2}$}\psfrag{P}{$c_{l-1}$}\psfrag{b}{$u_2$}
  \psfrag{h}{$u_1$}\psfrag{i}{$u_2$}\psfrag{j}{$u_{p-3}$}\psfrag{l}{$u_{p-2}$}
  \psfrag{f}{$u$}\psfrag{k}{$v$}\psfrag{r}{$w_2$}\psfrag{q}{$w_1$}\psfrag{t}{$w_{q-2}$}\psfrag{s}{$w_{q-3}$}
  \psfrag{A}{$a_1$}\psfrag{B}{$a_2$}\psfrag{C}{$a_{p-2}$}\psfrag{D}{$a_{p-1}$}\psfrag{E}{$b_1$}
  \psfrag{F}{$b_2$}\psfrag{G}{$b_{l-2}$}\psfrag{H}{$b_{l-1}$}\psfrag{8}{$v_i$}\psfrag{m}{$v_1$}
  \psfrag{n}{$v_2$}\psfrag{o}{$v_{l-3}$}\psfrag{p}{$v_{l-2}$}\psfrag{I}{$c_1$}\psfrag{J}{$c_2$}
  \psfrag{K}{$c_{q-2}$}\psfrag{L}{$c_{q-1}$}
  \includegraphics[width=140mm]{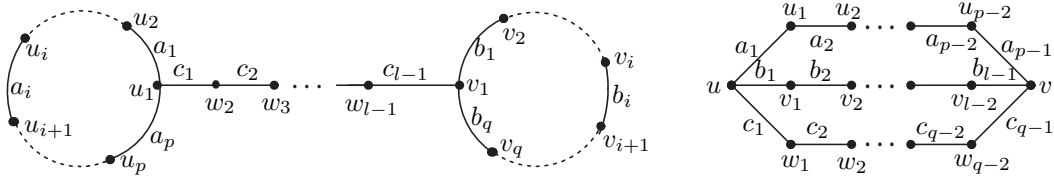}\\
  \caption{The weighted $\infty (p,l,q)$ graph and the weighted $\theta (p,l,q)$ graph.}
\end{center}
\end{figure}
\begin{figure}[h!]
\begin{center}
  \psfrag{2}{$u_1$}\psfrag{3}{$u_3$}\psfrag{1}{$u_2$}\psfrag{4}{$u_4$}
  \psfrag{a}{$a_4$}\psfrag{b}{$a_1$}\psfrag{c}{$a_2$}\psfrag{d}{$a_3$}
  \psfrag{e}{$c_1$}\psfrag{f}{$c_2$}\psfrag{g}{$c_k$}\psfrag{5}{$v_1$}
  \psfrag{6}{$v_2$}\psfrag{7}{$v_k$}\psfrag{G}{$G$}
  \includegraphics[width=50mm]{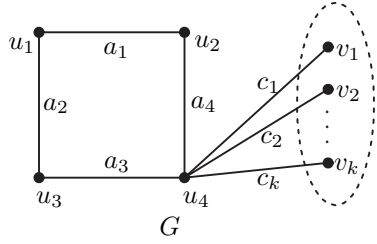}\\
  \caption{The weighted graph operation II of calculating positive and negative inertia index.}
\end{center}
\end{figure}

\begin{lem}\label{lem4.1}
Let $C_4=u_1u_3u_4u_2u_1$ be a weighted cycle as depicted in Fig. 3, where the number on each edge denotes the weight of the corresponding edge. Then
$$
i_+(C_4(u_4)\odot^kG)  =\left\{
  \begin{array}{ll}
    1+i_+(G+u_4), & \hbox{if $a_1a_3=a_2a_4;$} \\
    2+i_+(G), & \hbox{otherwise,}
  \end{array}
\right. i_-(C_4(u_4) \odot^kG)  = \left\{\begin{array}{ll}
    1+i_-(G+u_4), & \hbox{if $a_1a_3=a_2a_4;$} \\
    2+i_-(G), & \hbox{otherwise,}
  \end{array}\right.
$$
where $G+u$ is the induced weighted subgraph of $C_4(u_4)\odot ^kG.$
\end{lem}

\begin{proof}
Assume $N_{G}(u_4)=\{v_1,v_2,\ldots,v_k\}$ and $w(u_4v_i)=c_i> 0$ $(1\leq i\leq k)$; see Fig. 3. Then it follows that
\[A(C_4(u_4)\odot^kG)=
\left(\begin{array}{ccccccccccc}
0 & a_1& a_2 &0&0& 0&\ldots&0&0&\ldots & 0\\
a_1 & 0& 0 &a_4&0&0&\ldots& 0&0&\ldots & 0 \\
 a_2 & 0& 0 & a_3& 0&0&\ldots&0&0&\ldots & 0\\
0 & a_4& a_3 &0&c_1&c_2&\ldots& c_k&0&\ldots & 0 \\
0 & 0& 0 &c_1&0&a_{12}&\ldots& a_{1k}&a_{1(k+1)}&\ldots & a_{1(n-4)} \\
0 & 0& 0 &c_2&a_{12}&0&\ldots& a_{2k}&a_{2(k+1)}&\ldots & a_{2(n-4)} \\
\vdots & \vdots& \vdots &\vdots&\vdots&\vdots& \ddots&\vdots&\, &\vdots\\
0 & 0& 0 &c_k&a_{1k}&a_{2k}&\ldots&0&a_{k(k+1)}&\ldots & a_{k(n-4)}\\
0 & 0& 0 &0&a_{1(k+1)}&a_{2(k+1)}&\ldots&a_{k(k+1)}&0&\ldots & a_{(k+1)(n-4)}\\
\vdots & \vdots& \vdots & \vdots& \vdots&\vdots&\,&\vdots&\vdots&\ddots& \vdots\\
0& 0 &0&0&a_{1(n-4)}&a_{2(n-4)}&\ldots& a_{k(n-4)}& a_{(k+1)(n-4)}&\ldots & 0 \notag
\end{array}\right),\]
where the first $k+4$ rows and columns are labeled by $u_1,u_2,u_3,u_4,v_1,v_2,\ldots,v_k,$ respectively.

It is easy to show that $A(C_4(u_4)\odot^kG)$ is congruent to
\[
\left(\begin{array}{ccccccccccc}
0 & a_1& 0 &0&0& 0&\ldots&0&0&\ldots & 0\\
a_1 & 0& 0 &0&0&0&\ldots& 0&0&\ldots & 0 \\
 0 & 0& 0 & a_3-\frac{a_2a_4}{a_1}& 0&0&\ldots&0&0&\ldots & 0\\
0 & 0& a_3-\frac{a_2a_4}{a_1} &0&c_1&c_2&\ldots& c_k&0&\ldots & 0 \\
0 & 0& 0 &c_1&0&a_{12}&\ldots& a_{1k}&a_{1(k+1)}&\ldots & a_{1(n-4)} \\
0 & 0& 0 &c_2&a_{12}&0&\ldots& a_{2k}&a_{2(k+1)}&\ldots & a_{2(n-4)} \\
\vdots & \vdots& \vdots &\vdots&\vdots&\vdots& \ddots&\vdots&\, &\vdots\\
0 & 0& 0 &c_k&a_{1k}&a_{2k}&\ldots&0&a_{k(k+1)}&\ldots & a_{k(n-4)}\\
0 & 0& 0 &0&a_{1(k+1)}&a_{2(k+1)}&\ldots&a_{k(k+1)}&0&\ldots & a_{(k+1)(n-4)}\\
\vdots & \vdots& \vdots & \vdots& \vdots&\vdots&\,&\vdots&\vdots&\ddots& \vdots\\
0& 0 &0&0&a_{1(n-4)}&a_{2(n-4)}&\ldots& a_{k(n-4)}& a_{(k+1)(n-4)}&\ldots & 0 \notag
\end{array}\right).\]

If $a_1a_3=a_2a_4$, then $i_+(C_4(u_4)\odot^kG)  =
   1+i_+(G+u_4)$, $i_-(C_4(u_4)\odot^kG)  =
   1+i_-(G+u_4)$. Otherwise, $A(C_4(u_4)\odot^kG)$ is congruent to
\[
\left(\begin{array}{ccccccccccc}
0 & a_1& 0 &0&0& 0&\ldots&0&0&\ldots & 0\\
a_1 & 0& 0 &0&0&0&\ldots& 0&0&\ldots & 0 \\
 0 & 0& 0 & a_3-\frac{a_2a_4}{a_1}& 0&0&\ldots&0&0&\ldots & 0\\
0 & 0& a_3-\frac{a_2a_4}{a_1} &0&0&0&\ldots& 0&0&\ldots & 0 \\
0 & 0& 0 &0&0&a_{12}&\ldots& a_{1k}&a_{1(k+1)}&\ldots & a_{1(n-4)} \\
0 & 0& 0 &0&a_{12}&0&\ldots& a_{2k}&a_{2(k+1)}&\ldots & a_{2(n-4)} \\
\vdots & \vdots& \vdots &\vdots&\vdots&\vdots& \ddots&\vdots&\, &\vdots\\
0 & 0& 0 &0&a_{1k}&a_{2k}&\ldots&0&a_{k(k+1)}&\ldots & a_{k(n-4)}\\
0 & 0& 0 &0&a_{1(k+1)}&a_{2(k+1)}&\ldots&a_{k(k+1)}&0&\ldots & a_{(k+1)(n-4)}\\
\vdots & \vdots& \vdots & \vdots& \vdots&\vdots&\,&\vdots&\vdots&\ddots& \vdots\\
0& 0 &0&0&a_{1(n-4)}&a_{2(n-4)}&\ldots& a_{k(n-4)}& a_{(k+1)(n-4)}&\ldots & 0 \notag
\end{array}\right).\]
Then $i_+(C_4(u_4)\odot^kG)  =
   2+i_+(G)$, $i_-(C_4(u_4)\odot^kG)  =
   2+i_-(G)$.
\end{proof}

By the similar method, we can obtain the following lemma.
\begin{lem}\label{lem4.2}
Let $C_6$ be a weighted cycle of order 6 containing vertex $u$. Then
$i_+(C_6(u)\odot^kG) =3+i_+(G)$ and $i_-(C_6(u)\odot^kG) =3+i_-(G).$
\end{lem}

By Lemmas 3.4, 3.5, 3.10, 4.1 and 4.2, we can calculate the positive and negative inertia index of the weighted $\infty (p,l,q)$, if $p$ \,or $q\in \{4,6\}$. By the similar method, we can also determine the positive and negative inertia index of the weighted $\infty (p,l,q)$, if $p, q\in \{3,5\}$, as shown in Table 1, where the empty cell means there is no correlation between inertia index of $G$ and the weight set of $G$.
\begin{table}[h!]
  \centering
  \caption{Positive, negative inertia indices of graphs $\infty(p,l,q)\ (p,q\in \{3,5\}, 1\le l\le 5).$}\label{dd}\vspace{2mm}

\begin{tabular}{|c|c|c|c|}
  \hline
  Graphs $G$ & weighted conditions & $i_+(G)$ & $i_-(G)$  \\[3pt]\hline
  $\infty(3,1,3)$ &  & 2 & 3  \\[3pt]\hline
  \multirow{4}{16mm}{$\infty(3,2,3)$}
 &$4\frac{a_1b_1a_3b_3}{a_2b_2}-c_1^{2}>0$  & 2 & 4 \\[3pt]\cline{2-4}
 &$4\frac{a_1b_1a_3b_3}{a_2b_2}-c_1^{2}=0$  & 2 & 3 \\[3pt]\cline{2-4}
 &$4\frac{a_1b_1a_3b_3}{a_2b_2}-c_1^{2}<0$  & 3 & 3 \\[3pt]\hline
  $\infty(3,3,3)$ &  & 3 & 4  \\[3pt]\hline
  \multirow{4}{16mm}{$\infty(3,4,3)$}
 &$4\frac{a_1b_1a_3b_3}{a_2b_2}c_2^{2}-c_1^2c_3^2>0$  & 3 & 5 \\[3pt]\cline{2-4}
 &$4\frac{a_1b_1a_3b_3}{a_2b_2}c_2^{2}-c_1^2c_3^2=0$  & 3 & 4 \\[3pt]\cline{2-4}
 &$4\frac{a_1b_1a_3b_3}{a_2b_2}c_2^{2}-c_1^2c_3^2<0$  & 4 & 4 \\[3pt]\hline
 $\infty(3,5,3)$ &  & 4 & 5  \\[3pt]\hline
  \multirow{4}{16mm}{$\infty(3,1,5)$}
 &$\frac{a_1a_3}{a_2}-\frac{b_1b_3b_5}{b_2b_4}>0$  & 3 & 4 \\[3pt]\cline{2-4}
 &$\frac{a_1a_3}{a_2}-\frac{b_1b_3b_5}{b_2b_4}=0$  & 3 & 3 \\[3pt]\cline{2-4}
 &$\frac{a_1a_3}{a_2}-\frac{b_1b_3b_5}{b_2b_4}<0$  & 4 & 3 \\[3pt]\hline
 $\infty(3,2,5)$ &  & 4 & 4  \\[3pt]\hline
 \multirow{4}{16mm}{$\infty(3,3,5)$}
 &$\frac{a_1a_3}{a_2}c_2^2-\frac{b_1b_3b_5}{b_2b_4}c_1^2>0$  & 4 & 5 \\[3pt]\cline{2-4}
 &$\frac{a_1a_3}{a_2}c_2^2-\frac{b_1b_3b_5}{b_2b_4}c_1^2=0$  & 4 & 4 \\[3pt]\cline{2-4}
 &$\frac{a_1a_3}{a_2}c_2^2-\frac{b_1b_3b_5}{b_2b_4}c_1^2<0$  & 5 & 4 \\[3pt]\hline
 $\infty(3,4,5)$ &  & 5 & 5  \\[3pt]\hline
 \multirow{4}{16mm}{$\infty(3,5,5)$}
 &$\frac{a_1a_3}{a_2}c_2^2c_4^2-\frac{b_1b_3b_5}{b_2b_4}c_1^2c_3^2>0$  & 5 & 6 \\[3pt]\cline{2-4}
 &$\frac{a_1a_3}{a_2}c_2^2c_4^2-\frac{b_1b_3b_5}{b_2b_4}c_1^2c_3^2=0$  & 5 & 5 \\[3pt]\cline{2-4}
 &$\frac{a_1a_3}{a_2}c_2^2c_4^2-\frac{b_1b_3b_5}{b_2b_4}c_1^2c_3^2<0$  & 6 & 5 \\[3pt]\hline
  $\infty(5,1,5)$ &  & 5 & 4  \\[3pt]\hline
 \multirow{4}{16mm}{$\infty(5,2,5)$}
 &$4\frac{a_1b_1a_3b_3a_5b_5}{a_2b_2a_4b_4}-c_1^{2}>0$  & 6 & 4 \\[3pt]\cline{2-4}
 &$4\frac{a_1b_1a_3b_3a_5b_5}{a_2b_2a_4b_4}-c_1^{2}=0$  & 5 & 4 \\[3pt]\cline{2-4}
 &$4\frac{a_1b_1a_3b_3a_5b_5}{a_2b_2a_4b_4}-c_1^{2}<0$ & 5 & 5 \\[3pt]\hline
   $\infty(5,3,5)$ &  & 6 & 5  \\[3pt]\hline
 \multirow{4}{16mm}{$\infty(5,4,5)$}
 &$4\frac{a_1b_1a_3b_3a_5b_5}{a_2b_2a_4b_4}c_2^{2}-c_1^{2}c_3^2>0$  & 7 & 5 \\[3pt]\cline{2-4}
 &$4\frac{a_1b_1a_3b_3a_5b_5}{a_2b_2a_4b_4}c_2^{2}-c_1^{2}c_3^2=0$  & 6 & 5 \\[3pt]\cline{2-4}
 &$4\frac{a_1b_1a_3b_3a_5b_5}{a_2b_2a_4b_4}c_2^{2}-c_1^{2}c_3^2<0$ & 6 & 6 \\[3pt]\hline
 $\infty(5,5,5)$ &  & 7 & 6  \\\hline
\end{tabular}
\end{table}

From Lemmas 3.4, 3.5, 3.10, 4.1, 4.2 and Table 1, we can obtain the positive and negative inertia index of weighted $\infty(p,l,q)$, when $3\leq p\leq 6$, $3\leq q\leq 6$, $1\leq l\leq 5$. By Lemma 2.6, we obtain the following theorem immediately.
\begin{thm}\label{thm 4.3}
For each nonnegative integer $k, s$ and $t$, let $a = 4k+p$ $(3 \leq p \leq 6)$, $b = 4s+q$ ($3 \leq q \leq 6$), $c = 4t +l$
($1 \leq l \leq 5$). Then
$$
i_+(\infty(a, c, b)) = 2(k + s + t) + i_+(\infty(p, l, q)),\ \ i_-(\infty(a, c, b)) = 2(k + s + t) + i_-(\infty(p, l, q)).
$$
\end{thm}
\begin{lem}\label{lem4.6}
For weighted graph $\theta(3, 3, q),$ one has
$$
i_+(\theta(3, 3, q))  = \left\{
                          \begin{array}{ll}
                            i_+(C_{q+1}), & \hbox{if $a_1b_2=a_2b_1;$} \\
                            2+i_+(P_{q-2}), & \hbox{otherwise,}
                          \end{array}
                        \right.
i_-(\theta(3, 3, q))  = \left\{
                          \begin{array}{ll}
                            i_-(C_{q+1}), & \hbox{if $a_1b_2=a_2b_1;$} \\
                            2+i_-(P_{q-2}), & \hbox{otherwise.}
                          \end{array}
                        \right.
$$
\end{lem}
\begin{proof}
The weighted adjacency matrix $A(\theta(3, 3, q))$ has the following form
\[A(\theta(3, 3, q))=
\left(\begin{array}{cccccccc}
0 & x& a_1 &b_1&x_1& x_2&\ldots & x_{q-2}\\
x & 0& a_2 &b_2&y_1& y_2&\ldots & y_{q-2}\\
 a_1& a_2&0 &0&0&0&\ldots & 0\\
b_1 & b_2&0 &0&0&0&\ldots& 0 \\
x_1 & y_1&0 &0&0&a_{12}&\ldots&a_{1(q-2)} \\
x_2 & y_2&0 &0&a_{12}&0&\ldots&a_{2(q-2)} \\
\vdots & \vdots&\vdots &\vdots&\vdots&\vdots&\ddots&\vdots \\
x_{q-2} & y_{q-2}&0 &0&a_{1(q-2)}&a_{2(q-2)}&\ldots&0  \notag
\end{array}\right),\]
where $x>0$, if $q=2$; otherwise, $x=0$. And the first four rows and columns are labeled by $u,v,u_1,v_1,$ respectively.

It is easy to show that $A(\theta(3, 3, q))$ is congruent to
\[
\left(\begin{array}{cccccccc}
0 & x& a_1 &0&x_1& x_2&\ldots & x_{q-2}\\
x & 0& a_2 &b_2-\frac{b_1a_2}{a_1}&y_1& y_2&\ldots & y_{q-2}\\
 a_1& a_2&0 &0&0&0&\ldots & 0\\
0 & b_2-\frac{b_1a_2}{a_1}&0 &0&0&0&\ldots& 0 \\
x_1 & y_1&0 &0&0&a_{12}&\ldots&a_{1(q-2)} \\
x_2 & y_2&0 &0&a_{12}&0&\ldots&a_{2(q-2)} \\
\vdots & \vdots&\vdots &\vdots&\vdots&\vdots&\ddots&\vdots \\
x_{q-1} & y_{q-1}&0 &0&a_{1(q-2)}&a_{2(q-2)}&\ldots&0  \notag
\end{array}\right).\]

If $a_1b_2=a_2b_1$, $i_+(\theta(3, 3, q)) =i_+(\theta(3, 3, q)-v_1)=i_+(C_{q+1})$,  $i_-(\theta(3, 3, q)) =i_-(\theta(3, 3, q)-v_1)=i_+(C_{q+1})$. Otherwise, $A(\theta(3, 3, q))$ is congruent to
\[
\left(\begin{array}{cccccccc}
0 & 0& a_1 &0&0& 0&\ldots &0\\
0 & 0&0 &b_2-\frac{b_1a_2}{a_1}&0& 0&\ldots & 0\\
 a_1& 0&0 &0&0&0&\ldots & 0\\
0 & b_2-\frac{b_1a_2}{a_1}&0 &0&0&0&\ldots& 0 \\
0 & 0&0 &0&0&a_{12}&\ldots&a_{1(q-2)} \\
0 & 0&0 &0&a_{12}&0&\ldots&a_{2(q-2)} \\
\vdots & \vdots&\vdots &\vdots&\vdots&\vdots&\ddots&\vdots \\
0 & 0&0 &0&a_{1(q-2)}&a_{2(q-2)}&\ldots&0  \notag
\end{array}\right).\]
Therefore,
\begin{eqnarray*}
  i_+(\theta(3, 3, q)) &=& 2+i_+(\theta(3, 3,q)-u_1-v_1-u-v)=2+i_+(P_{q-2}), \\
  i_-(\theta(3, 3, q)) &=& 2+i_-(\theta(3, 3, q)-u_1-v_1-u-v)=2+i_-(P_{q-2}),
\end{eqnarray*}
as desired.
\end{proof}
By the similar argument, we can obtain the following lemmas.
\begin{lem}\label{lem4.7}
Let $G^{\ast}$ be a weighted cycle of order $q+2$ obtained from $\theta(4, 4, q)-v_1-v_2$ by replacing $a_3$ with $a_3+\frac{b_1a_2b_3}{a_1b_2}$. Then
$i_+(\theta(4, 4, q))=
   1+i_+(G^{\ast})$, $i_-(\theta(4, 4, q))=
   1+i_-(G^{\ast})$.
\end{lem}
\begin{lem}\label{lem4.8} For weighted graph $\theta(5, 5, q),$ one has
\begin{eqnarray*}
  i_+(\theta(5, 5, q)) &=& \left\{
                             \begin{array}{ll}
                               1+i_+(C_{q+3}), & \hbox{if $a_1b_2a_3b_4=b_1a_2b_3a_4;$} \\
                               2+i_+(P_{q+2}), & \hbox{otherwise.}
                             \end{array}
                           \right.
  \\
i_-(\theta(5, 5, q)) &=& \left\{
                             \begin{array}{ll}
                               1+i_-(C_{q+3}), & \hbox{if $a_1b_2a_3b_4=b_1a_2b_3a_4;$} \\
                               2+i_-(P_{q+2}), & \hbox{otherwise.}
                             \end{array}
                           \right.
\end{eqnarray*}
\end{lem}
\begin{lem}\label{lem4.9}
Let $G^{\ast}$ be a weighted cycle of order $q$ obtained from $\theta(2, 6, q)-v_1-v_2-v_3-v_4$ by replacing $a_1$ with $a_1+\frac{b_1b_3b_5}{b_2b_4}$. Then
$i_+(\theta(2, 6, q))=
   2+i_+(G^{\ast})$, $i_-(\theta(2, 6, q))=
   2+i_-(G^{\ast})$.
\end{lem}
\begin{lem}\label{lem4.10}For weighted graphs $\theta(2,4,3),\, \theta(2,4,5),\, \theta(2,3,5)$ and $\theta(3,4,5),$ one has
\begin{eqnarray*}
  i_+ (\theta(2,4,3)) &=& \left\{
                            \begin{array}{ll}
                              2, & \hbox{if $a_1b_2>b_1b_3;$} \\[5pt]
                              2, & \hbox{if $a_1b_2=b_1b_3;$}\\[5pt]
                              3,& \hbox{if $a_1b_2<b_1b_3,$}
                            \end{array}
                          \right.\ \
  i_- (\theta(2,4,3)) = \left\{
                            \begin{array}{ll}
                              3, & \hbox{if $a_1b_2>b_1b_3;$} \\[5pt]
                              2, & \hbox{if $a_1b_2=b_1b_3;$}\\[5pt]
                              2,& \hbox{if $a_1b_2<b_1b_3,$}
                            \end{array}
                          \right. \\
 i_+ (\theta(2,4,5)) &=& \left\{
                            \begin{array}{ll}
                              3, & \hbox{if $a_1b_2>b_1b_3;$} \\[5pt]
                              3, & \hbox{if $a_1b_2=b_1b_3;$}\\[5pt]
                              4,& \hbox{if $a_1b_2<b_1b_3,$}
                            \end{array}
                          \right.\ \
  i_- (\theta(2,4,5)) = \left\{
                            \begin{array}{ll}
                              4, & \hbox{if $a_1b_2>b_1b_3;$} \\[5pt]
                              3, & \hbox{if $a_1b_2=b_1b_3;$}\\[5pt]
                              3,& \hbox{if $a_1b_2<b_1b_3,$}
                            \end{array}
                          \right.  \\[5pt]
i_+(\theta(2,3,5))&=&i_-(\theta(2,3,5))=3,\,\,\,\,\, \ \ \ \ \ \ \, i_+(\theta(3,4,5))=i_-(\theta(3,4,5))=4.
\end{eqnarray*}
\end{lem}
From Lemmas 3.1, 3.4, 3.5 and 4.4-4.8, we can obtain the positive and negative inertia index of weighted $\infty(p,l,q)$, when $2\leq p\leq 6$, $2\leq q\leq 6$, $2\leq l\leq 6$. By Lemma 2.6, we obtain the following theorem immediately.
\begin{thm}\label{thm4.4}
For each nonnegative integer $k, s$ and $t$, let $a = 4k+p\, (2 \leq p \leq 5),\, b = 4s+q\ (2 \leq q \leq 5), c = 4t +l
(2 \leq l \leq 5).$ Then
$$
   i_+(\theta(a, c, b)) = 2(k + s + t) + i_+(\theta(p, l, q)),\ \ i_-(\theta(a, c, b)) = 2(k + s + t) + i_-(\theta(p, l, q)).
$$
\end{thm}
As we know, the connected bicyclic graphs can be partitioned into two classes: one is the set of all bipartite graphs each of
which contain an $\infty$-graph as an induced subgraph and the other is the set of all bipartite graphs each of which contain a $\theta$-graph as
an induced subgraph. The graph $\infty(p, l, q)$ (or $\theta(p, l, q)$) is called the \textit{base} of the corresponding bicyclic
graph $B$ which contains it, denoted by $\chi_{B}$.

Let $B$ be a weighted bicyclic graph with base $\chi_{B}$. For each vertex $v\in \chi_{B}$, let $B\{v\}$ be the induced
weighted connected subgraph of $B$ with the maximum possible number of vertices, which contains the vertex
$v$ and contains no other vertices of $\chi_{B}$. Then for all vertices $v\in \chi_{B}$, $B\{v\}$ is a weighted tree and $B$ is obtained by identifying the vertex $v$ of $B\{v\}$ with the vertex $v$ on $\chi_{B}$ . The bicyclic graph $B$ is called of Type I, if there
exists a vertex $v$ on $\chi_{B}$  such that $v$ is matched in $B\{v\}$; otherwise, $B$ is called of Type I\!I.

Similar to the discussion of weighted unicyclic graph, we have the following theorem.
\begin{thm}\label{thm4.1}
Let $B$ be a weighted bicyclic graph and let $\chi_{B}$ be the base of $B$.
\begin{wst}
\item[{\rm (i)}]
If $B$ is of type I and $v \in \chi_{B}$ is matched
in $B\{v\}$, then
$$
  i_+(B) =i_+(B\{v\}) + i_+(B - B\{v\}),\ \  i_-(B) =i_-(B\{v\}) + i_-(B - B\{v\}).
$$
And $B\{v\}$ is a tree, $B - B\{v\}$ is the union of weighted unicyclic graphs and trees.
\item[{\rm (ii)}]
If $B$ is of Type I\!I, then
$$
 i_+(B) = i_+(B - \chi_{B}) + i_+(\chi_{B}),\ \ i_-(B) = i_-(B- \chi_{B}) + i_-(\chi_{B}).
$$
\end{wst}
\end{thm}

\end{document}